\UseRawInputEncoding
\documentclass[12pt]{article}

\usepackage{mathrsfs}
\usepackage{amssymb}
\usepackage{amsmath}
\usepackage{amsfonts}
\usepackage{amstext}
\usepackage{amsthm}
\usepackage{graphicx}
\usepackage{subfig}
\usepackage{color}
\usepackage{indentfirst,latexsym,bm}
\usepackage{mathrsfs}
\usepackage{cite}
\usepackage[all]{xy}
\usepackage{newlfont}

\usepackage{enumerate}
\numberwithin{equation}{section}

\numberwithin{equation}{section}

\begin{document}

\title{Nonsymmetric sign-changing solutions for the partially overdetermined eigenvalue problem on a hollow cylinder 
\thanks{Research supported by NNSF of China (No. 12371110, 12301133).}}

\author{Guowei Dai
\thanks{School of Mathematical Sciences, Dalian University of Technology, Dalian, 116024, P.R. China
\newline
\text{~~~~ E-mail}: daiguowei@dlut.edu.cn}, Yingxin Sun\thanks{School of Mathematical Sciences, Dalian University of Technology, Dalian, 116024, P.R. China
\newline
\text{~~~~ E-mail}: sunyingxin2023@mail.dlut.edu.cn}, Yong Zhang \thanks{Corresponding author.
\newline
School of Mathematical Sciences, Jiangsu University, Zhenjiang, 212013, P.R. China
\newline
\text{~~~~ E-mail}: 18842629891@163.com} \\
}
\date{}
\maketitle

\renewcommand{\abstractname}{Abstract}

\begin{abstract}
{~
Fix $R\in(0,1)$ and $N\geq2$, and let $\lambda_k$ be the $k$th eigenvalue of the radial Dirichlet problem on the annulus
$A_R:=\{x\in\mathbb R^N:R<|x|<1\}$.  For every $k\geq1$ we construct a smooth local family of rotationally symmetric, axially periodic perturbations of the hollow cylinder $A_R\times\mathbb R$ for which
\begin{equation}
-\Delta u=\lambda u\quad\hbox{in }\Omega,\qquad
u=0\quad\hbox{on }\partial\Omega_{\mathrm{out}}\cup\partial\Omega_R,\qquad
\partial_\nu u=\hbox{constant}\quad\hbox{on }\partial\Omega_{\mathrm{out}}
\nonumber
\end{equation}
admits a nonradial periodic solution.  The branch bifurcates from the $k$th radial eigenfunction.  If $k\geq2$, the solution is sign-changing: it has precisely $k-1$ smooth nested nodal hypersurfaces and hence exactly $k$ nodal domains.  We also give an explicit connected-strip analogue when $N=1$.

The proof is based on a mean-free Dirichlet-to-Neumann operator.  A Sturm--Liouville Weyl function yields a unique simple crossing in the first spectral gap, so that the Crandall--Rabinowitz theorem produces a smooth local branch.  This gives a partially overdetermined construction on an unbounded hollow domain with disconnected complement; it lies outside the positivity and connected-complement hypotheses of the Berestycki--Caffarelli--Nirenberg conjecture.
}
\end{abstract}

\emph{Keywords:} Overdetermined problem; Bifurcation; Hollow cylinder; Bessel functions

{~\emph{AMS Subject Classification (2020):} 35B05; 35J25; 35P30; 58J55}

\tableofcontents

\section{Introduction}
\quad\, A domain $\Omega\subseteq\mathbb{R}^{N}$ is said to have the Pompeiu property if $f\equiv 0$ is the only continuous function
satisfying
\begin{equation}
\int_{\sigma(\Omega)}f(x)\,\text{d}x=0\nonumber
\end{equation}
for every rigid motion $\sigma$ of $\mathbb{R}^N$. The problem of classifying regions based on whether they have the Pompeiu property is called the Pompeiu problem.
The Pompeiu problem originated from harmonic analysis \cite{Pompeiu, Pompeiu1}.
Relevant conclusions can be applied to plasma physics \cite{Temam},
nuclear reactors \cite{Berenstein} and tomography \cite{Shepp, Smith}.
It is well known that a solid open ball of any radius $R > 0$ fails to have the Pompeiu property.
Conversely, it remains an open question whether a region without the Pompeiu property must necessarily be a ball.
In particular, Williams \cite{Williams} proposed the following conjecture.
\\ \\
\textbf{Williams conjecture.}
\emph{If $\partial\Omega$ is homeomorphic to the unit sphere in $\mathbb{R}^{N}$, then $\Omega$ has the Pompeiu
property if and only if it is not a ball.}\\

For a bounded regular domain $\Omega$, Williams conjecture holding is equivalent to the following Schiffer conjecture holding \cite{Williams}.
\\ \\
\textbf{Schiffer Conjecture.} \emph{Let $\Omega\subset \mathbb{R}^N$ be a bounded regular domain and assume that $u: \Omega \rightarrow \mathbb{R}$
is a solution to the problem
\begin{equation}
\left\{
\begin{array}{ll}
\Delta u+\lambda u=0\,\, &\text{in}\,\, \Omega,\\
u=1 &\text{on}\,\, \partial \Omega,\\
\partial_\nu u=0 &\text{on}\,\, \partial \Omega,
\end{array}
\right.\nonumber
\end{equation}
where $\nu$ is the unit outer normal vector on $\partial \Omega$. Then $\Omega$ is a ball and $u$ is radially symmetric.}\\

The Pompeiu problem or Schiffer conjecture is one of Yau's famous list of problems \cite[Problem 80]{Yau}. As far as we know, no one gave a thoroughly positive answer to this problem on bounded domain up to now. However, very recently, Fall, Minlend and Weth \cite{FallMW} provided some counterexamples to Schiffer conjecture on unbounded domains by using the method of bifurcation. Enciso, Fern\'{a}ndez, Ruiz and Sicbaldi \cite{EncisoFRS}  provided other counterexamples to Schiffer conjecture on  bounded planar doubly connected domains.
In addition, Berenstein \cite{Berenstein} also proposed the following variation of Schiffer conjecture by exchanging boundary conditions.
\\ \\
\textbf{Berenstein conjecture.}
\emph{Let $\Omega$ be a bounded $C^{2,\alpha}$ domain in $\mathbb{R}^{N}$ with $\alpha\in(0,1)$. If there exists a nontrivial
solution $u$ of the overdetermined eigenvalue problem
\begin{equation}\label{overdeterminedeigenvalueproblem}
\left\{
\begin{array}{ll}
\Delta u+\lambda u=0\,\, &\text{in}\,\, \Omega,\\
u=0 &\text{on}\,\, \partial \Omega,\\
\partial_\nu u=\text{const} &\text{on}\,\, \partial \Omega,
\end{array}
\right.
\end{equation}
then $\Omega$ is a ball, where $\nu$ is the unit outer normal on $\partial \Omega$.}\\

For $N=2$, Berenstein \cite{Berenstein} proved that the existence of infinitely many eigenvalues for (\ref{overdeterminedeigenvalueproblem}) is equivalent to $\Omega$ being a disk. Similar result was also deduced
to the Poincar\'{e} disk with the hyperbolic metric by Berenstein and Yang \cite{Berenstein1}. In \cite{Berenstein2}, Berenstein and Yang further showed that Berenstein's result was valid for any dimension.

Berenstein conjecture is also closely related to the result
of Serrin \cite{Serrin}. {~More precisely, if $\Omega$ is a bounded connected smooth domain and there exists
$u>0$ satisfying $\Delta u+1=0$ in $\Omega$, $u=0$ on $\partial\Omega$, and
$\partial_\nu u=\mathrm{const}$ on $\partial\Omega$, then $\Omega$ is a ball.}
Since Serrin's famous work \cite{Serrin}, overdetermined elliptic problem has attracted a lot of attention from many mathematicians. While the most investigation of overdetermined boundary value problems was essentially limited to proofs
that solutions need to be radial in cases that could be handled using the method of
moving planes \cite{Gidas}, for instance \cite{Aftalion, BCN, ChenLi, Gidas, Pucci, Reichel, RRS} (the lists are just a very small sample and far from being complete). In particular, Berestycki, Caffarelli and Nirenberg \cite{BCN}
consider the following more general overdetermined elliptic problem
\begin{equation}\label{ffunction}
\left\{
\begin{array}{ll}
\Delta u+f(u)=0\,\, &\text{in}\,\, \Omega,\\
u=0 &\text{on}\,\, \partial \Omega,\\
\partial_\nu u=\text{const.} &\text{on}\,\, \partial \Omega,
\end{array}
\right.
\end{equation}
where $f$ is a given locally Lipschitz function.
Under some suitable assumptions on $f$ and domain $\Omega$, Berestycki, Caffarelli and Nirenberg proved that if problem (\ref{ffunction}) admits a smooth and bounded solution, then $\Omega\subseteq \mathbb{R}^N$ is a half-space. Meanwhile, they \cite{BCN} also proposed the following famous conjecture.
\\ \\
\textbf{BCN conjecture.}
If $f$ is a Lipschitz function on a
domain $\Omega$ in $\mathbb{R}^N$ such that $\mathbb{R}^N\setminus \overline{\Omega}$ is connected, then the existence of a bounded positive solution
to (\ref{ffunction}) implies that $\Omega$ is either a ball, a half-space, a generalized cylinder
$B^m\times \mathbb{R}^{N-m}$ ($B^m$ is a round ball in $\mathbb{R}^m$), or the complement of one of them.\\

It is well known that many confirmed answers to the BCN conjecture have been given in bounded domain. In addition,
if $\Omega$ is $C^3$ and uniformly Lipschitz epigraph of $\mathbb{R}^2$ or $\mathbb{R}^3$, Farina and Valdinoci \cite{Farina} proved that
there is no solution $u\in C^2\left(\overline{\Omega}\right) \cap L^\infty (\Omega)$ of
\begin{equation}\label{positivesoution}
\left\{
\begin{array}{ll}
\Delta u+\lambda u=0\,\, &\text{in}\,\, \Omega,\\
u>0 &\text{in}\,\, \Omega,\\
u=0 &\text{on}\,\, \partial \Omega,\\
\partial_\nu u=\text{const} &\text{on}\,\, \partial \Omega.
\end{array}
\right.
\end{equation}
That is to say, the necessary condition for problem (\ref{positivesoution}) possessing a solution is that $\Omega$ is half-space, which can be seen as a confirmed answer to BCN or Berenstein conjecture on unbounded domain. For certain types of functions $f(u)$, BCN conjecture is also true for exterior domains (see \cite{Aftalion, Reichel}). However, the first counterexample to BCN conjecture on unbounded domain was constructed by Sicbaldi \cite{Sicbaldi} via showing that the cylinder $B_1\times \mathbb{R}$ with $N\geq2$ can be perturbed to an unbounded domain whose boundary is a periodic
hypersurface of revolution with respect to the $\mathbb{R}$-axis, where $B_1$ is the unit ball of $\mathbb{R}^N$ centered on the origin.
After that, Schlenk and Sicbaldi \cite{Schlenk} proved that the above conclusion is also valid for $N=1$ and these extremal domains belong to a smooth bifurcation family of domains. Nontrivial epigraph domain emanating from the half-space has been obtained in \cite{Del} for $N\geq8$. Ros, Ruiz and Sicbaldi \cite{Ros1} constructed a counter-example to BCN conjecture on exterior domain which is the perturbation of the complement of a ball.

Note that these mentioned counterexamples mainly involve a positive solution, and $\mathbb{R}^N\setminus \overline{\Omega}$
is connected. Recently, Minlend \cite{Minlend}, as well as Dai and Zhang \cite{DaiZ}, found sign-changing solutions on unbounded bifurcating domains, while Ruiz \cite{Ruiz} and Wheeler \cite{Wheeler} obtained sign-changing solutions on bounded bifurcating domains. We also refer the reader to \cite{LPS,PRS} for bifurcating domains arising from an overdetermined eigenvalue problem in cylinders. {~A natural related question is whether a nontrivial unbounded hollow domain, whose complement is disconnected, can support a nonsymmetric sign-changing solution when the Neumann overdetermination is imposed only on the outer boundary.  This is a partially overdetermined analogue of the problems above, but it is not a counterexample to the BCN conjecture: both positivity and connectedness of the complement are hypotheses of that conjecture.  The main goal of this paper is to answer this partially overdetermined question affirmatively.}

{~
\medskip
\noindent\textbf{Relation to earlier work and novelty.}
The geometry and the boundary data considered here differ from the existing
periodic constructions in three essential respects.  Dai and Zhang
\cite{DaiZ} perturb a solid cylinder and impose the full overdetermined
condition on its entire boundary.  Minlend \cite{MinlendPeriodic} considers
periodic domains with two boundary components on which prescribed Neumann
values of opposite signs are imposed.  By contrast, our inner cylinder is
fixed and carries only the Dirichlet condition; the constant Neumann datum is
required solely on the moving outer component.  The recent construction of
Lian, Pacella and Sicbaldi \cite{LPS} concerns positive eigenfunctions in
bounded hypographs contained in a half-cylinder, whereas our domains are
unbounded periodic hollow cylinders and the branches issued from
$\lambda_k$, $k\geq2$, are sign-changing.

The main analytical contribution is not merely the replacement of a ball by
an annulus.  The nonlinear normal trace is projected onto a fixed mean-zero
codomain, which removes the missing constant equation without changing the
overdetermined problem.  Moreover, a scalar Sturm--Liouville Weyl function
proves that, for every radial level and every $N\geq2$, the first Fourier mode
has exactly one zero in the first spectral gap, that this zero is simple, and
that no higher Fourier mode belongs to the kernel there.  This yields a
one-dimensional kernel and a genuine Crandall--Rabinowitz branch, together
with the exact persistence of all $k-1$ nodal hypersurfaces.  These features
are specific to the partially overdetermined hollow geometry and are not
consequences of the cited cylinder results.
}

More precisely, we consider domains of the form $\Omega:=\Omega_{out}\setminus\overline{\Omega}_{R}$, where
$$\overline{\Omega}_{R}=\left\{(x,t)\in \mathbb{R}^N\times \mathbb{R}: |x|\leq R\right\}$$
is a fixed solid cylinder and $\Omega_{\mathrm{out}}$ is a smooth perturbation of the unit cylinder $\Omega_{1}=\left\{(x,t)\in \mathbb{R}^N\times \mathbb{R}: |x|\leq 1\right\}$ with $\overline{\Omega}_{R}\subset \Omega_{\mathrm{out}}$.
Thus, $\Omega$ is a nontrivial hollow cylindrical domain in $\mathbb{R}^{N+1}$ (see Figure 1). We study the following partially overdetermined eigenvalue problem:
\begin{equation}\label{Hollowcylindereigenvalueproblem}
\left\{
\begin{array}{ll}
\Delta u+\lambda u=0\,\, &\text{in}\,\, \Omega,\\
u=0 &\text{on}\,\, \partial \Omega_{R}\cup \partial \Omega_{out},\\
\partial_\nu u=\text{const} &\text{on}\,\, \partial \Omega_{out},
\end{array}
\right.
\end{equation}
For any $k\in \mathbb{N}$, let $\lambda_k$ be the $k$th eigenvalue of the Dirichlet Laplacian restricted to radial functions on $A_R:=\left\{x\in \mathbb{R}^N:0<R<\vert x\vert<1\right\}$.
{~For every fixed radial level $k$, we obtain a smooth local family of unbounded domains bifurcating from
the hollow cylinder $A_R\times \mathbb{R}$ such that problem (\ref{Hollowcylindereigenvalueproblem}) has a nonsymmetric solution. More precisely, we establish the following results.}
\\ \\
{~
\textbf{Theorem 1.1. $(N\geq 2)$}
\itshape Fix $R\in(0,1)$, $N\geq2$, $k\geq1$, and $\alpha\in(0,1)$.  Let
$\mathcal{C}^{2,\alpha}_{\mathrm{even},0}(\mathbb{R}/2\pi\mathbb{Z})$
denote the even $2\pi$-periodic $\mathcal{C}^{2,\alpha}$ functions of mean zero.
There exist $\varepsilon>0$, a number $T_{1,*}=T_{1,*}(N,R,k)$ satisfying
\begin{equation}
 T_{1,*}\in\left(\frac{2\pi}{\sqrt{\lambda_k}},+\infty\right)
 \quad\text{if }k=1,
 \qquad
 T_{1,*}\in\left(\frac{2\pi}{\sqrt{\lambda_k}},
 \frac{2\pi}{\sqrt{\lambda_k-\lambda_1}}\right)
 \quad\text{if }k\geq2,
 \nonumber
\end{equation}
and smooth maps
\[
s\longmapsto (w_s,T_s,\lambda_s)
\quad\text{from }(-\varepsilon,\varepsilon)
\text{ into }
\mathcal{C}^{2,\alpha}_{\mathrm{even},0}(\mathbb{R}/2\pi\mathbb{Z})
\times(0,\infty)\times(0,\infty),
\]
such that $w_0=0$, $T_0=T_{1,*}$, $\lambda_0=\lambda_k$, and
$w_s$ is $L^2$-orthogonal to $\cos t$.  Set
\[
v_s(\theta)=s\cos\theta+s w_s(\theta)
\]
and
\[
\Omega_s=\left\{(x,t)\in\mathbb R^N\times\mathbb R:
R<|x|<1+v_s\!\left(\frac{2\pi t}{T_s}\right)\right\}.
\]
Then problem (\ref{Hollowcylindereigenvalueproblem}), with $\lambda=\lambda_s$,
admits a $T_s$-periodic solution
$u_s\in\mathcal C^{2,\alpha}(\overline{\Omega_s})$.
For $s\neq0$ sufficiently small, $\Omega_s$ is a nontrivial perturbation of
$A_R\times\mathbb R$ and $u_s$ is nonradial.

If $k=1$, $u_s$ can be chosen positive.  If $k\geq2$, there exist precisely
$k-1$ smooth $T_s$-periodic functions $r_{j,s}$ satisfying
\[
R<r_{1,s}(t)<\cdots<r_{k-1,s}(t)
<1+v_s\!\left(\frac{2\pi t}{T_s}\right),
\qquad
u_s(r_{j,s}(t),t)=0.
\]
These graphs are all the interior nodal hypersurfaces of $u_s$; consequently
$u_s$ has exactly $k$ nodal domains.
}
\begin{figure}[ht]
\centering
\includegraphics[width=0.7\textwidth]{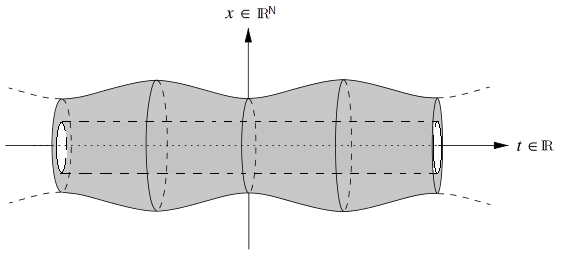}
\caption{The bifurcation domain of the hollow cylinder for $N\geq 2$.}
\end{figure}
~\\

In the special case $N=1$, we are able to obtain a stronger result as follows.
~\\
~\\
{~
\textbf{Theorem 1.2. $(N=1)$}
\itshape Fix $R\in(0,1)$, $k\geq1$, and $\alpha\in(0,1)$, and put
\[
T_{1,*}=\frac{4(1-R)}{\sqrt{4k^2-1}}.
\]
There exist $\varepsilon>0$ and smooth maps
\[
s\longmapsto (w_s,T_s,\lambda_s)
\]
defined for $|s|<\varepsilon$, with $w_s\in
\mathcal C^{2,\alpha}_{\mathrm{even},0}(\mathbb R/2\pi\mathbb Z)$,
$w_0=0$, $T_0=T_{1,*}$,
$\lambda_0=k^2\pi^2/(1-R)^2$, and
$\int_0^{2\pi}w_s(\theta)\cos\theta\,d\theta=0$, such that, for
\[
v_s(\theta)=s\cos\theta+s w_s(\theta),
\qquad
\Omega_s^+=\left\{(x,t)\in\mathbb R^2:
R<x<1+v_s\!\left(\frac{2\pi t}{T_s}\right)\right\},
\]
the corresponding partially overdetermined problem has a nontrivial
$T_s$-periodic solution $u_s\in
\mathcal C^{2,\alpha}(\overline{\Omega_s^+})$ with eigenvalue $\lambda_s$.
It is positive for $k=1$; for $k\geq2$ it has precisely $k-1$ periodic
nodal curves and exactly $k$ nodal domains.  Reflecting $\Omega_s^+$ and
$u_s$ across $x=0$ gives the two-component symmetric formulation
$R<|x|<1+v_s(2\pi t/T_s)$.
}\\

{~To prove Theorem 1.1 we choose the $k$th radial Dirichlet
eigenfunction on $A_R$ as the trivial state and formulate the outer
overdetermination by subtracting the mean of its Dirichlet-to-Neumann
trace.  The resulting map takes values in a fixed mean-zero space.  A
Sturm--Liouville Weyl function identifies a unique simple characteristic
value in the first spectral gap and verifies the transversality condition.
The Crandall--Rabinowitz theorem then gives the smooth branch stated above.
For $k\geq2$ its eigenfunction is sign-changing, and the complement of the
hollow domain is disconnected; these facts place the construction outside
the hypotheses of the BCN conjecture.}
We also would like to mention that some independent interests on partially overdetermined elliptic boundary value problems are considered in \cite{Farina1,Fragala,Fragala1}. In fact, the problem (\ref{Hollowcylindereigenvalueproblem}) considered in this paper is also regarded as a partially overdetermined elliptic boundary value problem because the two boundary conditions are put on outer boundary, however one just boundary condition is put on inner boundary. These problems appear quite naturally in many different phenomena in Physics, see \cite{Sirakov} for more details.

{~
The annular cross-section introduces a one-dimensional radial
Sturm--Liouville problem.  Although its eigenfunctions can be written in
terms of two Bessel functions, no delicate comparison of individual
Bessel zeros is needed for bifurcation.  The decisive quantity is the
logarithmic derivative of the normalized radial solution; its derivative
has a definite sign by a Lagrange identity.

The paper is organized as follows.  Section 2 gives the nonlinear
mean-free formulation, computes its linearization, and records the Bessel
representation of the radial eigenfunctions.  Section 3 proves existence
and uniqueness of the simple crossing.  Section 4 applies the
Crandall--Rabinowitz theorem and proves stability of the nodal
hypersurfaces.  Section 5 treats the connected one-dimensional strip and
computes the bifurcation period explicitly.
}

\section{Preliminaries}

\quad\,  We first show an elementary result on zero-Dirichlet Laplacian eigenvalue problems on cylinder in the following.

{~
\subsection{The radial eigenvalue problem on the cylinder}

For $N\geq2$ let
\[
A_R=\{x\in\mathbb R^N:R<|x|<1\}.
\]
Restriction of the Dirichlet Laplacian to $O(N)$-invariant functions gives
the regular Sturm--Liouville problem
\begin{equation}\label{radial-sl}
 -(r^{N-1}\phi')'=\lambda r^{N-1}\phi,\qquad
 \phi(R)=\phi(1)=0.
\end{equation}
Its eigenvalues are simple and form a strictly increasing sequence
\[
0<\lambda_1<\lambda_2<\cdots\longrightarrow+\infty.
\]
We normalize the corresponding eigenfunction $\phi_k$ by
\begin{equation}\label{normalization}
 \int_{A_R}\phi_k^2\,dx=\frac1{2\pi},
 \qquad \phi_k'(1)<0,
 \qquad b_k:=-\phi_k'(1)>0.
\end{equation}
By the Sturm oscillation theorem, $\phi_k$ has exactly $k-1$ simple
zeros
\[
R<\rho_1<\cdots<\rho_{k-1}<1.
\]
The equation at $r=1$ also gives the identity
\begin{equation}\label{phi-second}
 \phi_k''(1)=-(N-1)\phi_k'(1)=(N-1)b_k.
\end{equation}

For $T>0$, the function $(x,t)\mapsto\phi_k(|x|)$ is a
$T$-periodic eigenfunction on
\[
C_{R,1}^T=A_R\times(\mathbb R/T\mathbb Z).
\]
After the change of variable $\theta=2\pi t/T$, its $L^2$ norm on
$A_R\times\mathbb S^1$ is one.  The connected $N=1$ problem is recorded
separately in Section 5.
}

{~
\subsection{A mean-free nonlinear formulation}

Write $\mathbb S^1=\mathbb R/2\pi\mathbb Z$ and set
\begin{align*}
X&=\left\{v\in\mathcal C^{2,\alpha}(\mathbb S^1):
v(-\theta)=v(\theta),\ \int_0^{2\pi}v(\theta)\,d\theta=0\right\},\\
Y&=\left\{q\in\mathcal C^{1,\alpha}(\mathbb S^1):
q(-\theta)=q(\theta),\ \int_0^{2\pi}q(\theta)\,d\theta=0\right\}.
\end{align*}
Define the end of the first spectral gap by
\begin{equation}\label{gap-end}
 \widehat T_k=
 \begin{cases}
 +\infty,&k=1,\\[2mm]
 \displaystyle\frac{2\pi}{\sqrt{\lambda_k-\lambda_1}},&k\geq2.
 \end{cases}
\end{equation}
For $T\in(0,\widehat T_k)$ the eigenvalue $\lambda_k$ is simple in the
space of functions which are $O(N)$-invariant in $x$ and even in $t$.
Indeed, the spectrum in this symmetry class at the straight cylinder is
\[
 \lambda_j+\left(\frac{2m\pi}{T}\right)^2,
 \qquad j\geq1,\quad m\geq0,
\]
and the definition of $\widehat T_k$ excludes every equality with
$\lambda_k$ when $m\geq1$.

For $v\in X$ sufficiently small, let
\[
C_{R,1+v}^T=
 \left\{(x,t)\in\mathbb R^N\times(\mathbb R/T\mathbb Z):
 R<|x|<1+v\!\left(\frac{2\pi t}{T}\right)\right\}.
\]
A radial cut-off diffeomorphism maps this domain to $C_{R,1}^T$ and fixes
the inner boundary.  The implicit function theorem for the pulled-back
Dirichlet operator and the normalization
\[
 \int_{C_{R,1+v}^T}\phi_{k,v,T}^2\,dx\,dt=\frac{T}{2\pi}
\]
give, near every $(0,T_0)$ with $T_0\in(0,\widehat T_k)$, a unique smooth
eigenpair $(\phi_{k,v,T},\lambda_{k,v,T})$ in the above symmetry class,
with the sign fixed by \eqref{normalization}.  At $v=0$ this pair is
$(\phi_k,\lambda_k)$.

Let $\nu_{v,T}$ be the outer unit normal and pull the outer normal
derivative back to $\mathbb S^1$:
\[
 \mathcal D(v,T)(\theta)=
 \left.\partial_{\nu_{v,T}}\phi_{k,v,T}(x,t)\right|_{
 |x|=1+v(\theta),\,t=T\theta/(2\pi)}.
\]
With
\[
 Pq=q-\frac1{2\pi}\int_0^{2\pi}q(\eta)\,d\eta,
\]
define
\begin{equation}\label{definition-F}
 F(v,T)=P\mathcal D(v,T).
\end{equation}
Then $F$ is a smooth map from a neighborhood of
$\{0\}\times(0,\widehat T_k)$ in $X\times\mathbb R$ to the fixed space
$Y$, and $F(0,T)=0$.  Most importantly,
\begin{equation}\label{F-equivalence}
 F(v,T)=0
 \quad\Longleftrightarrow\quad
 \partial_{\nu_{v,T}}\phi_{k,v,T}\ \hbox{is constant on the outer boundary}.
\end{equation}
Subtracting the actual mean in \eqref{definition-F}, rather than a fixed
normal derivative, is what makes both the nonlinear equation and its
linearization have the correct codomain.

{~
\medskip
\noindent\textbf{Lemma 2.1 (smooth eigenpair branch and smoothness of
$F$).}
\emph{Fix $T_0\in(0,\widehat T_k)$.  There are neighborhoods
$\mathcal U\subset X$ of $0$ and $I\subset(0,\widehat T_k)$ of $T_0$
such that the normalized eigenpair
\[
 (v,T)\longmapsto(\phi_{k,v,T},\lambda_{k,v,T})
\]
is uniquely defined and smooth on $\mathcal U\times I$.  After pullback
to the fixed cylinder, the eigenfunction depends smoothly in
$\mathcal C^{2,\alpha}$ and the eigenvalue depends smoothly in
$\mathbb R$.  Consequently the map
$F:\mathcal U\times I\to Y$ in \eqref{definition-F} is smooth.}

\medskip
\noindent\emph{Proof.}
Put $D=A_R\times\mathbb S^1$ and choose
$\chi\in\mathcal C^\infty([R,1])$ satisfying
$\chi=0$ near $R$ and $\chi=1$ near $1$.  For $\|v\|_{\mathcal
C^{2,\alpha}}$ small, the map
\begin{equation}\label{fixed-domain-map}
 \Psi_v(r\omega,\theta)
 =\big((r+\chi(r)v(\theta))\omega,\theta\big),
 \qquad \omega\in\mathbb S^{N-1},
\end{equation}
is a $\mathcal C^{2,\alpha}$ diffeomorphism from $D$ onto the
$\theta$-rescaled periodic cell.  It fixes $r=R$ and maps $r=1$ onto
$r=1+v(\theta)$.  In these coordinates the Euclidean metric is pulled
back to
\begin{equation}\label{pulled-metric}
 g_{v,T}=\Psi_v^*\left(
 dr^2+r^2g_{\mathbb S^{N-1}}
 +\left(\frac{T}{2\pi}\right)^2d\theta^2\right).
\end{equation}
Let
\[
 E=\{U\in\mathcal C^{2,\alpha}(\overline D):
 U|_{\partial D}=0,\ U\text{ is }O(N)\text{-invariant and even in }\theta\}
\]
and let $Z$ be the analogous subspace of
$\mathcal C^{0,\alpha}(\overline D)$.  Define
\[
 L(v,T)U=-\Delta_{g_{v,T}}U,\qquad
 d\mu_{v,T}=\frac{2\pi}{T}\,d\operatorname{vol}_{g_{v,T}}.
\]
The coefficients of $L(v,T)$ are smooth functions of
$v,v',v''$ and $T$ on the open set where \eqref{fixed-domain-map} is
a diffeomorphism.  The usual multiplication and inversion properties
of Hölder spaces therefore show that
\[
 (v,T,U)\longmapsto L(v,T)U
 \quad\text{is smooth from }\mathcal U\times I\times E\text{ to }Z.
\]

Consider
\begin{equation}\label{EIFT-map}
 \mathcal G(v,T,U,\lambda)=
 \left(L(v,T)U-\lambda U,\
 \int_DU^2\,d\mu_{v,T}-1\right).
\end{equation}
At $(0,T_0,\phi_k,\lambda_k)$ this map vanishes.  Its derivative in
$(U,\lambda)$ is
\begin{equation}\label{EIFT-derivative}
 (\xi,\mu)\longmapsto
 \left((L(0,T_0)-\lambda_k)\xi-\mu\phi_k,\
 2\langle\phi_k,\xi\rangle_{L^2(d\mu_{0,T_0})}\right).
\end{equation}
This is an isomorphism from $E\times\mathbb R$ onto
$Z\times\mathbb R$.  Indeed, on the orthogonal complement of
$\phi_k$ the first component is invertible by the Fredholm alternative
and the Dirichlet Schauder estimate, while the $\phi_k$ component and
the last scalar equation determine $\mu$ and the component of $\xi$
along $\phi_k$.

For clarity, the required simplicity follows directly from the first
spectral gap.  The spectrum of $L(0,T_0)$ in the chosen symmetry class
is
\[
 \lambda_j+\left(\frac{2m\pi}{T_0}\right)^2,\qquad
 j\geq1,\quad m\geq0.
\]
If $m\geq1$, then for $k\geq2$
\[
 \lambda_j+\left(\frac{2m\pi}{T_0}\right)^2
 \geq\lambda_1+\left(\frac{2\pi}{T_0}\right)^2
 >\lambda_k,
\]
and the same conclusion is immediate for $k=1$.  Hence the only
eigenfunction at $\lambda_k$ in this symmetry class is $\phi_k$.

The Banach-space implicit function theorem applied to
\eqref{EIFT-map} now gives unique smooth maps
$U_{k,v,T}\in E$ and $\lambda_{k,v,T}\in\mathbb R$.  The sign is fixed
by requiring
$\langle U_{k,v,T},\phi_k\rangle_{L^2(d\mu_{0,T_0})}>0$.
Pushing $U_{k,v,T}$ forward by \eqref{fixed-domain-map} gives
$\phi_{k,v,T}$ and the normalization used above.

Finally, after pullback to $r=1$, the unit normal and the normal trace
are smooth algebraic functions of $v,v'$ and the first derivatives of
$U_{k,v,T}$.  The trace theorem in Hölder spaces thus yields a smooth
map
\[
 (v,T)\longmapsto\mathcal D(v,T)
 \quad\text{from }\mathcal U\times I
 \text{ to }\mathcal C^{1,\alpha}_{\mathrm{even}}(\mathbb S^1).
\]
Since $P$ is a fixed bounded projection onto $Y$,
$F=P\mathcal D$ is smooth. \hfill$\square$
}

\subsection{Linearization and Fourier multipliers}

Let $w\in X$.  Its shape derivative is described by the unique solution
$\psi_w$ with zero $\theta$-average of
\begin{equation}\label{shape-derivative}
\left\{
\begin{aligned}
 &\left(\Delta_x+\left(\frac{2\pi}{T}\right)^2\partial_\theta^2
       +\lambda_k\right)\psi_w=0
       &&\text{in }A_R\times\mathbb S^1,\\
 &\psi_w=0 &&\text{on }r=R,\\
 &\psi_w=-\phi_k'(1)w=b_kw &&\text{on }r=1.
\end{aligned}\right.
\end{equation}
The first variation of the eigenvalue is zero because $w$ has mean zero.
The standard Hadamard calculation, followed by the projection $P$, gives
\begin{equation}\label{linearization}
 H_Tw:=D_vF(0,T)w
 =P\left(\partial_r\psi_w(1,\cdot)+\phi_k''(1)w\right).
\end{equation}

{~
\medskip
\noindent\textbf{Proposition 2.2 (shape derivative and the mean-free
linearization).}
\emph{For every $T\in(0,\widehat T_k)$ and $w\in X$, the first
variation of the eigenvalue in the direction $w$ is zero.  The
Eulerian derivative of the normalized eigenfunction is the unique
solution $\psi_w$ of \eqref{shape-derivative} with zero
$\theta$-average, and \eqref{linearization} holds.  In fact the
expression inside $P$ already has mean zero.}

\medskip
\noindent\emph{Proof.}
Let $\Gamma_1=\{r=1\}\times\mathbb S^1$ and deform only $\Gamma_1$
with velocity field $V=w(\theta)e_r$; the velocity vanishes near
$r=R$.  Write $(\phi_s,\lambda_s)$ for the normalized eigenpair on
the domain corresponding to $v=sw$, and denote its Eulerian
derivative at $s=0$ by $\psi$.

The differentiated weak eigenvalue identity, followed by integration by
parts, gives the Hadamard formula
\begin{equation}\label{Hadamard-lambda}
 \dot\lambda=
 -\frac{\displaystyle\int_{\Gamma_1}
       (\partial_\nu\phi_k)^2w\,dS_T}
       {\displaystyle\int_{C_{R,1}^T}\phi_k^2\,dV_T}.
\end{equation}
Here $dS_T$ and $dV_T$ are the surface and volume measures in the
$T$-periodic cell.  Since $\partial_\nu\phi_k=\phi_k'(1)=-b_k$ is
constant on $\Gamma_1$ and the remaining factors in $dS_T$ are
independent of $\theta$, the defining condition
$\int_0^{2\pi}w(\theta)\,d\theta=0$ implies
\begin{equation}\label{lambda-dot-zero}
 \dot\lambda=0.
\end{equation}

Differentiation of
$\Delta\phi_s+\lambda_s\phi_s=0$ at fixed interior points, together
with \eqref{lambda-dot-zero}, yields
\[
 \left(\Delta_x+\left(\frac{2\pi}{T}\right)^2\partial_\theta^2
       +\lambda_k\right)\psi=0.
\]
The inner boundary is fixed, so $\psi=0$ at $r=R$.  Differentiating
$\phi_s(1+sw(\theta),\theta)=0$ on the outer boundary gives
\[
 \psi(1,\theta)+\phi_k'(1)w(\theta)=0,
\]
which is exactly the outer condition in
\eqref{shape-derivative}.

There remains a possible multiple of $\phi_k$ in the zero Fourier
mode.  The derivative of the normalization has no boundary term,
because $\phi_k=0$ on both boundary components, and therefore
\[
 \int_{C_{R,1}^T}\phi_k\psi\,dV_T=0.
\]
Taking the $\theta$-average of the equation for $\psi$ shows that this
average is a Dirichlet solution at the eigenvalue $\lambda_k$, hence
it is a multiple of $\phi_k$.  The last orthogonality identity forces
that multiple to vanish.  Thus $\psi=\psi_w$ is precisely the unique
zero-average solution of \eqref{shape-derivative}.

It remains to differentiate the outer normal trace.  At $s=0$ the
variation of the unit normal is tangential:
\[
 \dot\nu=-\nabla_{\Gamma_1}w.
\]
Since $\phi_k$ is radial,
$\nabla_{\Gamma_1}\phi_k=0$, so the term
$\nabla\phi_k\cdot\dot\nu$ vanishes.  Differentiating both the
eigenfunction and the evaluation point $r=1+sw(\theta)$ gives
\begin{equation}\label{normal-trace-variation}
 \left.\frac{d}{ds}\right|_{s=0}\mathcal D(sw,T)
 =\partial_r\psi_w(1,\theta)+\phi_k''(1)w(\theta).
\end{equation}
Applying the fixed projection $P$ proves
\eqref{linearization}.  Finally,
\[
 \int_0^{2\pi}\partial_r\psi_w(1,\theta)\,d\theta
 =\partial_r\int_0^{2\pi}\psi_w(1,\theta)\,d\theta=0,
\]
and $\int_0^{2\pi}w=0$.  Hence the right-hand side of
\eqref{normal-trace-variation} already belongs to $Y$.
\hfill$\square$
}

For
\[
 w(\theta)=\sum_{m\geq1}a_m\cos(m\theta),
 \qquad
 \alpha_m(T)=\lambda_k-\left(\frac{2m\pi}{T}\right)^2,
\]
write
\[
 \psi_w(r,\theta)=\sum_{m\geq1}a_mc_m(r;T)\cos(m\theta).
\]
Then
\begin{equation}\label{cm-problem}
 c_m''+\frac{N-1}{r}c_m'+\alpha_m(T)c_m=0,\qquad
 c_m(R;T)=0,\quad c_m(1;T)=b_k,
\end{equation}
and
\begin{equation}\label{sigma-def}
 H_T(\cos m\theta)=\sigma_m(T)\cos m\theta,\qquad
 \sigma_m(T)=c_m'(1;T)+\phi_k''(1).
\end{equation}
For $T\in(0,\widehat T_k)$ one has $\alpha_m(T)<\lambda_1$, so
\eqref{cm-problem} is uniquely solvable.

For completeness, this dependence is analytic.  Put
$\ell(r)=b_k(r-R)/(1-R)$ and $c_m=\ell+u$.  The function $u$ belongs to
the linear Dirichlet space
$\mathcal C^{2,\alpha}_0([R,1])$ and satisfies an equation
$(A-\alpha_m(T))u=f_{\alpha_m(T)}$, where
\[
 Au=-r^{1-N}(r^{N-1}u')'.
\]
Since $\alpha_m(T)<\lambda_1$, the resolvent is bounded.  The resolvent
identity (or its locally convergent Neumann series on this linear space)
shows that $(A-\alpha)^{-1}$ is analytic for
$\alpha\notin\{\lambda_j:j\geq1\}$.  Consequently $c_m'(1;T)$ and
$\sigma_m(T)$ are analytic on the stated interval.
}

\subsection{Some conclusions on Bessel functions}

\quad\, For $\tau\in \mathbb{R}$ and $s\in \mathbb{R}$, the Bessel function of the first kind is defined by
\begin{equation}
J_\tau(s)=\sum_{m=0}^\infty\frac{(-1)^m\left(\frac{s}{2}\right)^{2m+\tau}}{m!\Gamma(\tau+m+1)},\nonumber
\end{equation}
which is the solution of the differential equation
\begin{equation}
s^2f''(s)+sf'(s)+\left(s^2-\tau^2\right)f(s)=0.\nonumber
\end{equation}
The Bessel function of the second kind is defined by
\begin{equation}
Y_\tau(s)=\left\{
\begin{array}{ll}
\frac{J_\tau(s)\cos(\tau \pi)-J_{-\tau}(s)}{\sin(\tau \pi)}\,\, &\text{if}\,\, \tau\not\in \mathbb{Z},\\
\lim_{p\rightarrow \tau}\frac{J_p(s)\cos(p \pi)-J_{-p}(s)}{\sin(p \pi)} &\text{if}\,\, \tau\in \mathbb{Z}.
\end{array}
\right.\nonumber
\end{equation}
In particular, for any integer $n$, one has that
\begin{equation}
Y_n(s)=\frac{1}{\pi}\frac{\partial J_\tau(s)}{\partial \tau}\Big|_{\tau=n}+\frac{(-1)^n}{\pi}\frac{\partial J_\tau(s)}{\partial \tau}\Big|_{\tau=-n}.\nonumber
\end{equation}
Thus the general solution of the Bessel's differential equation is given by
\begin{equation}
AJ_\tau(s)+BY_\tau(s)\nonumber
\end{equation}
for any constants $A,B\in \mathbb{R}$.

The Bessel functions $J_\tau(s)$ and $Y_\tau(s)$ have the following differentiation formulas
\begin{equation}\label{weifengongshi3}
\frac{d}{d s}\left(\frac{J_\tau(s)}{s^{\tau}}\right)=-\frac{J_{\tau+1}(s)}{s^{\tau}}
\end{equation}
and
\begin{equation}\label{weifengongshi4}
\frac{d}{d s}\left(\frac{Y_\tau(s)}{s^{\tau}}\right)=-\frac{Y_{\tau+1}(s)}{s^{\tau}}.
\end{equation}
It follows that
\begin{equation}\label{weifengongshi5}
J_{\tau+1}(s)=\frac{\tau}{s}J_{\tau}(s)-J_{\tau}'(s)
\end{equation}
and
\begin{equation}\label{weifengongshi6}
Y_{\tau+1}(s)=\frac{\tau}{s}Y_{\tau}(s)-Y_{\tau}'(s)
\end{equation}
for $s\neq0$.
If $\tau$ is real, then each of $J_\tau(s)$ and $Y_\tau(s)$ has an infinite number of positive real zeros. All of
these zeros are simple, the $m$-th positive zero of $J_\tau(s)$ and $Y_\tau(s)$ are denoted by $j_{\tau,m}$ and $y_{\tau,m}$, respectively.
When $\tau\geq0$, the zeros interlace according to the
inequalities
{~
\begin{equation}\label{interlaceproperty2}
0<y_{\tau,1}<j_{\tau,1}<y_{\tau,2}<j_{\tau,2}<\cdots .
\end{equation}
}
When $\tau\geq-1$ the zeros of $J_\tau(s)$ are all real.

For $\tau\in \mathbb{R}$, the modified Bessel function of the first kind is defined by
\begin{equation}
I_\tau(s)=\sum_{m=0}^\infty\frac{\left(\frac{s}{2}\right)^{2m+\tau}}{m!\Gamma(\tau+m+1)},\nonumber
\end{equation}
which is the solution of the following modified Bessel's differential equation
\begin{equation}
s^2f''(s)+sf'(s)-\left(s^2+\tau^2\right)f(s)=0.\nonumber
\end{equation}
The modified Bessel function of the second kind is
{~
\begin{equation}
K_\tau(s)=\left\{
\begin{array}{ll}
\displaystyle\frac{\pi}{2}\frac{I_{-\tau}(s)-I_\tau(s)}
{\sin(\pi\tau)} &\text{if }\tau\notin\mathbb Z,\\[3mm]
\displaystyle\lim_{p\to\tau}\frac{\pi}{2}
\frac{I_{-p}(s)-I_p(s)}{\sin(\pi p)}
&\text{if }\tau\in\mathbb Z.
\end{array}
\right.\nonumber
\end{equation}
}
which is also the solution of the modified Bessel's differential equation.

For any $\tau\in \mathbb{R}$, we have that
\begin{equation}\label{weifengongshi}
\frac{d}{d s}\left(\frac{I_\tau(s)}{s^{\tau}}\right)=\frac{I_{\tau+1}(s)}{s^{\tau}}
\end{equation}
and
\begin{equation}\label{weifengongshi1}
\frac{d}{d s}\left(\frac{K_\tau(s)}{s^{\tau}}\right)=-\frac{K_{\tau+1}(s)}{s^{\tau}}.
\end{equation}
The function $I_\tau(s)$ and $K_\tau(s)$ have the following asymptotic properties
\begin{equation}\label{asymptoticproperties}
I_\tau(s)\sim \frac{e^s}{\sqrt{2\pi s}},\,\,s\rightarrow+\infty
\end{equation}
and
\begin{equation}\label{asymptoticproperties1}
K_\tau(s)\sim \sqrt{\frac{\pi}{2s} }e^{-s},\,\,s\rightarrow+\infty.
\end{equation}
The function $I_\tau(s)$ and $K_\tau(s)$ are linearly independent.
The general solution of the modified Bessel's differential equation is given by
\begin{equation}
AI_\tau(s)+BK_\tau(s)\nonumber
\end{equation}
for any constants $A,B\in \mathbb{R}$.

{~
\subsection{A nonsingular Bessel representation of $\phi_k$}

The following formula is useful for identifying the radial spectrum, even
though the crossing argument in Section 3 will use only
Sturm--Liouville theory.  Put
\[
 \tau=\frac{N-2}{2},\qquad q_k=\sqrt{\lambda_k}.
\]
After the substitution $s=q_kr$ and
$\widehat\phi(s)=s^\tau\phi(s/q_k)$,
equation \eqref{radial-sl} becomes Bessel's equation of order
$\tau$.  Hence, for a nonzero normalizing constant $C_k$,
\begin{equation}\label{bessel-eigenfunction-correct}
 \phi_k(r)=C_kr^{-\tau}
 \left[
   Y_\tau(q_kR)J_\tau(q_kr)
   -J_\tau(q_kR)Y_\tau(q_kr)
 \right].
\end{equation}
This determinant combination automatically vanishes at $r=R$ and does
not divide by a possibly vanishing Bessel value.  The outer Dirichlet
condition is exactly
\begin{equation}\label{annular-characteristic}
 Y_\tau(q_kR)J_\tau(q_k)
 -J_\tau(q_kR)Y_\tau(q_k)=0.
\end{equation}
Conversely, every positive root of \eqref{annular-characteristic} gives
a radial Dirichlet eigenvalue.  Simplicity and ordering are consequences
of the Sturm--Liouville problem \eqref{radial-sl}, so no additional sign
assumption on the individual Bessel factors is required.
}

{~
\section{The Weyl function and the unique simple crossing}
\label{sec:weyl}

For $\alpha<\lambda_1$, let $y_\alpha$ be the solution of
\begin{equation}\label{weyl-ivp}
 y_\alpha''+\frac{N-1}{r}y_\alpha'+\alpha y_\alpha=0,\qquad
 y_\alpha(R)=0,\quad y_\alpha'(R)=1,
\end{equation}
and define its Weyl logarithmic derivative by
\begin{equation}\label{m-definition}
 m(\alpha)=\frac{y_\alpha'(1)}{y_\alpha(1)}.
\end{equation}
Sturm comparison shows that $y_\alpha(1)>0$ for
$\alpha<\lambda_1$, so \eqref{m-definition} has no pole on this
interval.

\medskip
\noindent\textbf{Proposition 3.1.}
\emph{The function $m$ is analytic and strictly decreasing on
$(-\infty,\lambda_1)$.  More precisely,
\begin{equation}\label{m-derivative}
 m'(\alpha)=-
 \frac{\displaystyle\int_R^1r^{N-1}y_\alpha(r)^2\,dr}
      {y_\alpha(1)^2}<0.
\end{equation}
Moreover,
\[
 m(0)>0,\qquad \lim_{\alpha\uparrow\lambda_1}m(\alpha)=-\infty.
\]
Consequently there is a unique $\alpha_*\in(0,\lambda_1)$ such that
\begin{equation}\label{alpha-star}
 m(\alpha_*)=-(N-1).
\end{equation}}

\medskip
\noindent\emph{Proof.}
Analytic dependence follows from the initial-value problem.  Set
$z_\alpha=\partial_\alpha y_\alpha$.  Differentiating
\eqref{weyl-ivp}, multiplying the equations for $z_\alpha$ and
$y_\alpha$ by $y_\alpha$ and $z_\alpha$, respectively, and subtracting
gives
\[
 \left[r^{N-1}(z_\alpha'y_\alpha-z_\alpha y_\alpha')\right]'
 =-r^{N-1}y_\alpha^2.
\]
The initial data of $z_\alpha$ vanish at $R$; evaluation at $1$ yields
\eqref{m-derivative}.  When $\alpha=0$, the solution of
\eqref{weyl-ivp} is positive and increasing on $(R,1]$, hence
$m(0)>0$.

At $\alpha=\lambda_1$, $y_{\lambda_1}$ is a positive multiple of the
first radial eigenfunction with $y_{\lambda_1}'(1)<0$.  The same
Lagrange identity at $\lambda_1$ gives
\[
 \partial_\alpha y_\alpha(1)\big|_{\alpha=\lambda_1}\,
 y_{\lambda_1}'(1)
 =\int_R^1r^{N-1}y_{\lambda_1}^2\,dr>0.
\]
Thus $\partial_\alpha y_\alpha(1)|_{\lambda_1}<0$, and the quotient
in \eqref{m-definition} tends to $-\infty$ from the left.  Strict
monotonicity now proves existence and uniqueness in
\eqref{alpha-star}. \hfill$\square$

\medskip
Using \eqref{cm-problem} and \eqref{phi-second}, we obtain the exact
multiplier formula
\begin{equation}\label{sigma-m-weyl}
 \sigma_m(T)=b_k\left[m\!\left(
 \lambda_k-\left(\frac{2m\pi}{T}\right)^2\right)+(N-1)\right].
\end{equation}
Define
\begin{equation}\label{T-star-def}
 T_{1,*}=\frac{2\pi}{\sqrt{\lambda_k-\alpha_*}}.
\end{equation}
Since $0<\alpha_*<\lambda_1$, this number satisfies the bounds stated
in Theorem 1.1.  Formula \eqref{alpha-star} gives
$\sigma_1(T_{1,*})=0$.  If $m\geq2$, then
\[
 \lambda_k-\left(\frac{2m\pi}{T_{1,*}}\right)^2
 <\lambda_k-\left(\frac{2\pi}{T_{1,*}}\right)^2=\alpha_*.
\]
Together with \eqref{m-derivative}, this implies
\begin{equation}\label{higher-positive}
 \sigma_m(T_{1,*})>0\qquad(m\geq2).
\end{equation}
The crossing is transversal because
\begin{equation}\label{sigma-transverse}
 \sigma_1'(T_{1,*})
 =b_km'(\alpha_*)\frac{8\pi^2}{T_{1,*}^3}<0.
\end{equation}
We have therefore proved both uniqueness of the first-mode crossing in
the first spectral gap and simplicity of the full even mean-zero kernel.

{~
\medskip
\noindent\textbf{Lemma 3.2 (high-frequency symbol and Fredholm
property).}
\emph{Let $J\Subset(0,\widehat T_k)$.  Uniformly for $T\in J$,
\begin{equation}\label{symbol-asymptotic}
 \sigma_m(T)=b_k\frac{2m\pi}{T}+O(1)\qquad(m\to\infty).
\end{equation}
More precisely, if
\[
 a_T=\frac{2\pi b_k}{T},\qquad
 \rho_m(T)=\sigma_m(T)-a_Tm,
\]
then, for every integer $\ell\geq0$,
\begin{equation}\label{symbol-differences}
 |\Delta_m^\ell\rho_m(T)|
 \leq C_{\ell,J}(1+m)^{-\ell}.
\end{equation}
Consequently $H_T:X\to Y$ is Fredholm of index zero.  It is
self-adjoint with respect to the $L^2(\mathbb S^1)$ pairing, and at
$T=T_{1,*}$,
\begin{equation}\label{kernel-range}
 \ker H_{T_{1,*}}=\operatorname{span}\{\cos\theta\},\qquad
 \operatorname{Ran}H_{T_{1,*}}
 =\left\{q\in Y:\int_0^{2\pi}q(\theta)\cos\theta\,d\theta=0\right\}.
\end{equation}}

\medskip
\noindent\emph{Proof.}
For all sufficiently large $m$, set
\[
 \kappa_m(T)=
 \sqrt{\left(\frac{2m\pi}{T}\right)^2-\lambda_k}.
\]
Then \eqref{cm-problem} becomes
\[
 c_m''+\frac{N-1}{r}c_m'-\kappa_m^2c_m=0.
\]
Let $p=(N-1)/2$ and write $c_m=r^{-p}Z$.  The function $Z$ solves
\begin{equation}\label{Z-equation}
 Z''-\big(\kappa_m^2+Q(r)\big)Z=0,\qquad
 Q(r)=\frac{(N-1)(N-3)}{4r^2}.
\end{equation}
Let $Z_\kappa$ denote the solution of \eqref{Z-equation} with
$Z_\kappa(R)=0$ and $Z_\kappa'(R)=1$.  Variation of constants gives
the Volterra equation
\begin{equation}\label{Volterra-Z}
 Z_\kappa(r)=\frac{\sinh(\kappa(r-R))}{\kappa}
 +\int_R^r\frac{\sinh(\kappa(r-s))}{\kappa}
 Q(s)Z_\kappa(s)\,ds.
\end{equation}
Since $R>0$, the potential $Q$ and all its derivatives are bounded on
$[R,1]$.  Iteration of \eqref{Volterra-Z}, followed by differentiation
in $r$, gives
\begin{equation}\label{log-derivative-asymptotic}
 \frac{Z_\kappa'(1)}{Z_\kappa(1)}
 =\kappa\coth(\kappa(1-R))+O(\kappa^{-1})
 =\kappa+O(\kappa^{-1})
\end{equation}
as $\kappa\to+\infty$.  The constants are uniform for $T\in J$.
Because logarithmic derivatives are unaffected by normalization and
$c_m(1)=b_k$, we have
\[
 c_m'(1;T)=b_k\left(
 \frac{Z_{\kappa_m}'(1)}{Z_{\kappa_m}(1)}-p\right).
\]
Together with
\[
 \kappa_m(T)=\frac{2m\pi}{T}+O(m^{-1}),
\]
this proves \eqref{symbol-asymptotic}.

For completeness, differentiating \eqref{Volterra-Z} repeatedly with
respect to $\kappa$ gives, after the same Volterra estimate,
\[
 \left|\partial_\kappa^\ell\left(
 \frac{Z_\kappa'(1)}{Z_\kappa(1)}-\kappa\right)\right|
 \leq C_{\ell,J}(1+\kappa)^{-\ell},
 \qquad \ell\geq0.
\]
Since $\kappa_m(T)$ is a classical symbol of order one in $m$, the
finite-difference mean value formula yields
\eqref{symbol-differences}.

Let $|D|$ be the periodic multiplier
$|D|\cos(m\theta)=m\cos(m\theta)$.  Equations
\eqref{symbol-differences} give the decomposition
\begin{equation}\label{HT-decomposition}
 H_T=a_T|D|+K_T,
\end{equation}
where $K_T$ is a zero-order periodic Fourier multiplier.  The
standard dyadic proof of the periodic multiplier theorem, using
\eqref{symbol-differences} and summation by parts on each dyadic
block, shows that
\[
 K_T:\mathcal C^{2,\alpha}(\mathbb S^1)
 \longrightarrow\mathcal C^{2,\alpha}(\mathbb S^1)
\]
is bounded.  Hence $K_T:X\to Y$ is compact by the compact embedding
$\mathcal C^{2,\alpha}\hookrightarrow\mathcal C^{1,\alpha}$.
On the even mean-zero spaces, $|D|:X\to Y$ is an isomorphism; its
inverse has multiplier $1/m$.  Thus \eqref{HT-decomposition} proves
that $H_T$ is Fredholm and
\[
 \operatorname{ind}H_T=\operatorname{ind}(a_T|D|)=0.
\]

All multipliers $\sigma_m(T)$ are real, so the Fourier representation
also proves self-adjointness.  At $T=T_{1,*}$,
\eqref{higher-positive} and $\sigma_1(T_{1,*})=0$ give the stated
kernel.  If $q\in Y$ has vanishing $\cos\theta$ coefficient, define
\[
 w=\sum_{m\geq2}\frac{q_m}{\sigma_m(T_{1,*})}\cos(m\theta).
\]
The reciprocal symbol is of order $-1$, by
\eqref{symbol-differences}, so $w\in X$ and
$H_{T_{1,*}}w=q$.  Conversely every element of the range is
$L^2$-orthogonal to the kernel by self-adjointness.  This proves
\eqref{kernel-range}. \hfill$\square$
}

\section{Bifurcation and nodal stability}

\noindent\textbf{Proof of Theorem 1.1.}
The map $F:X\times(0,\widehat T_k)\to Y$ is smooth,
$F(0,T)=0$, and its derivative in $v$ is $H_T$.
Equations \eqref{kernel-range} and \eqref{sigma-transverse} imply
\[
 D_{Tv}F(0,T_{1,*})[\cos\theta]
 =\sigma_1'(T_{1,*})\cos\theta
 \notin\operatorname{Ran}H_{T_{1,*}}.
\]
All hypotheses of the Crandall--Rabinowitz theorem \cite{Crandall} are
therefore satisfied.  There are $\varepsilon>0$ and smooth maps
$s\mapsto(w_s,T_s)$ such that
\[
 w_0=0,\qquad T_0=T_{1,*},\qquad
 \int_0^{2\pi}w_s(\theta)\cos\theta\,d\theta=0,
\]
and, with $v_s=s(\cos\theta+w_s)$,
\[
 F(v_s,T_s)=0\qquad(|s|<\varepsilon).
\]
The eigenvalue construction in Section 2 also gives the smooth function
$\lambda_s=\lambda_{k,v_s,T_s}$ with $\lambda_0=\lambda_k$.
By \eqref{F-equivalence}, the corresponding eigenfunction satisfies all
three boundary conditions in
\eqref{Hollowcylindereigenvalueproblem}.  Periodic lifting from the cell
produces the solution on $\Omega_s$.  Since
$v_s/s\to\cos\theta$ in $\mathcal C^{2,\alpha}$, the boundary is
nonconstant and the solution is nonradial for $s\ne0$ sufficiently
small.

It remains to justify the assertion about signs and nodal sets.  Pull
the bifurcating eigenfunction back to
$A_R\times\mathbb S^1$ by the same cut-off diffeomorphism used in
Section 2 and denote it by $U_s$.  Smooth dependence gives
\begin{equation}\label{eigenfunction-convergence}
 U_s\longrightarrow\phi_k
 \quad\text{in }\mathcal C^{2,\alpha}
 (\overline{A_R\times\mathbb S^1}).
\end{equation}
Each interior zero $\rho_j$ of $\phi_k$ is simple.  The implicit
function theorem and \eqref{eigenfunction-convergence} therefore give
one and only one periodic graph $r=r_{j,s}(\theta)$ near $\rho_j$.
On the complement of small disjoint neighborhoods of these graphs,
$\phi_k$ is bounded away from zero; the nonzero normal derivatives at
$r=R,1$ control the boundary collars.  Hence
\eqref{eigenfunction-convergence} rules out every additional zero.
The $k-1$ nested graphs separate the periodic hollow cylinder into
exactly $k$ connected nodal regions.  When $k=1$, the branch is the
continuation of the simple principal eigenpair on the periodic cell, so
its eigenfunction can be chosen positive.  This completes the proof.
\hfill$\square$

\section{The connected one-dimensional strip}

Let $N=1$ and work on the connected component $R<x<1$; put
$L=1-R$.  With the same normalization and orientation as in
\eqref{normalization},
\begin{equation}\label{one-d-eigenpair}
 \lambda_k=\frac{k^2\pi^2}{L^2},\qquad
 \phi_k(x)=\frac{(-1)^{k+1}}{\sqrt{\pi L}}
 \sin\!\left(\frac{k\pi(x-R)}{L}\right),
 \qquad b_k=-\phi_k'(1)>0.
\end{equation}
The factor in \eqref{one-d-eigenpair} makes the $L^2$ norm on
$(R,1)\times\mathbb S^1$ equal to one.

The nonlinear construction and the linearization of Section 2 apply
without change, now with $N-1=0$.  For $0<\alpha<\lambda_1$ the solution
of \eqref{weyl-ivp} is
\[
 y_\alpha(x)=\frac{\sin(\sqrt\alpha(x-R))}{\sqrt\alpha},
 \qquad
 m(\alpha)=\sqrt\alpha\cot(L\sqrt\alpha).
\]
Its unique zero in $(0,\lambda_1)$ is
\[
 \alpha_*=\frac{\pi^2}{4L^2}.
\]
Substitution in \eqref{T-star-def} gives
\begin{equation}\label{one-d-period}
 T_{1,*}
 =\frac{2\pi}{\sqrt{\lambda_k-\alpha_*}}
 =\frac{4(1-R)}{\sqrt{4k^2-1}}.
\end{equation}
Formula \eqref{m-derivative} still gives
$\sigma_1'(T_{1,*})\ne0$, while strict monotonicity gives
$\sigma_m(T_{1,*})>0$ for every $m\geq2$.  The
Crandall--Rabinowitz and nodal-stability arguments of Section 4 prove
Theorem 1.2 on the connected strip.  Reflection by $x\mapsto-x$
preserves the equation and the two boundary conditions on each
component, yielding the optional two-component formulation.
}

{\bf Data Availability}: This article has no associated data.
~\\
~\\
{\bf Declarations}
~\\
~\\
{\bf Conflict of interest}: The authors declare that there are no conflict of interest.
\bibliographystyle{amsplain}
\makeatletter
\renewcommand{\@biblabel}[1]{[#1]}
\makeatother


\providecommand{\bysame}{\leavevmode\hbox to3em{\hrulefill}\thinspace}
\providecommand{\MR}{\relax\ifhmode\unskip\space\fi MR }
\providecommand{\MRhref}[2]{%
  \href{http://www.ams.org/mathscinet-getitem?mr=#1}{#2}
}
\providecommand{\href}[2]{#2}

\end{document}